\documentclass[10pt]{article}
\usepackage{graphicx}
\usepackage{amsmath,amssymb,amsthm,amsfonts}
\usepackage{amssymb}

\newtheorem{thm}{Theorem}[section]

\newtheorem{rem}[thm]{Remark}
\newcommand{\bremark}{\begin{rem} \textup}
\newcommand{\eremark}{\end{rem} }

\newcommand{\C}{{\mathbb C}}
\newcommand{\cuad}{{\sqcap\kern-.68em\sqcup}}

\newcommand{\R}{{\mathbb{R}}}
\newcommand{\N}{{\mathbb{N}}}

\def\ni{\noindent}
\def\proof{{\ni \bf \underline{Proof:} }}
\def\qed{{\unskip\nobreak\hfil\penalty50
         \hskip2em\hbox{}\nobreak\hfil\mbox{\rule{1ex}{1ex} \qquad}
           \parfillskip=0pt
           \finalhyphendemerits=0\par }}

\begin{document}

\title{Linearized theory for entire solutions of a singular Liouvillle equation}

\author{Manuel del Pino \footnote{Departamento de Ingenier\'ia Matem\'atica and CMM, Universidad de Chile, Casilla 170, Correo 3, Santiago, Chile. E-mail: delpino@dim.uchile.cl. Author supported by grants Fondecyt 1070389 and FONDAP (Chile)}\quad
Pierpaolo Esposito \footnote{Dipartimento di Matematica,
Universit\`a degli Studi ``Roma Tre", Largo S. Leonardo Murialdo, 1
-- 00146 Roma, Italy. E-mail: esposito@mat.uniroma3.it. Author partially supported by by FIRB-IDEAS (2008), project ``Geometrical aspects in PDEs".}\quad
Monica Musso \footnote{Departamento de Matematica, Pontificia Universidad Catolica de
Chile, Avenida Vicuna Mackenna 4860, Macul, Santiago, Chile. E-mail: mmusso@mat.puc.cl. Author  supported by Fondecyt grant 1040936 (Chile), and by M.U.R.S.T., project
``Metodi variazionali e topologici nello studio di fenomeni non lineari''.}}

\maketitle

\begin{abstract}  \noindent We discuss invertibility properties for entire finite-energy solutions of the regularized version of a singular Liouvillle equation.
\end{abstract}



\vskip 0.2truein

\section{The linear theory}
In this short note, we want to show the following result
\begin{thm} \label{linop}
Let $c \in \C$ and $N\in \N$. The kernel of the operator
$$L: \: \phi\: \to\: L(\phi):=\Delta
\phi+\frac{8(N+1)^2|z|^{2N}}{(1+|z^{N+1}-c|^2)^2}\phi$$ in
$L^\infty(\R^2)$ has the simple form
$$\{\phi \in L^\infty(\R^2):\:L(\phi)=0\}=\hbox{Span }\{\frac{1-|z^{N+1}-c|^2}{1+|z^{N+1}-c|^2},
\frac{\hbox{Re }(z^{N+1}-c)}{1+|z^{N+1}-c|^2}, \frac{\hbox{Im
}(z^{N+1}-c)}{1+|z^{N+1}-c|^2}\},$$ where $z^{N+1}$ denotes the
complex $(N+1)-$power.
\end{thm}
The case $N=0$ was already known (see \cite{BaPa}), based on a Fourier expansion of the function $\phi$. The aim here is to adapt the argument to the more difficult case $N\in \mathbb{N}$.\\
\proof Let us recall the Liouville formula: given a holomorphic
function $f$ on $\C$, the function
\begin{equation}\label{Liouville}
\ln \frac{8 | f'(z)|^2}{(1+|f(z)|^2)^2}
\end{equation}
solves the equation $\Delta U+e^U=0$ in the set $\{z \in \C \ /\
f'(z)\not=0 \}$. If now $f'$ has a zero at the origin of
multiplicity $N$, the function
\begin{equation}\label{Liouvillebis}
\ln \frac{8 |f'(z)|^2}{(1+|f(z)|^2)^2}-\ln |z|^{2N}
\end{equation}
solves the equation $\Delta U+|z|^{2N}e^U=0$ in the set $\{z \in
\C \setminus\{0\} \ /\ f'(z)\not=0 \}$. The choice
$f(z)=z^{N+1}(1+\tau z^k)-c$, $k \geq 0$, leads to a family
$$U_{\tau,k}(z) =  \ln \frac{8(N+1)^2|1+\tau \frac{N+1+k}{N+1}z^k|^2 }
{(1+|z^{N+1}(1+\tau z^k)-c|^2)^2}, \:\:\tau \in \C,$$
of solutions of $\Delta U+|z|^{2N}e^U=0$ in $\C \setminus \{z \in
\C: 1+\tau \frac{N+1+k}{N+1}z^k=0\}$. The derivative of
$U_{\tau,k}$ in $\tau$ at $\tau=0$:
$$\phi_k:=\partial_\tau U_{\tau,k} \Big|_{\tau=0}=z^k
\left(\frac{N+1+k}{N+1}-2\frac{z^{N+1}\overline{z^{N+1}-c}}{1+
|z^{N+1}-c|^2} \right)$$ solves $L(\phi_k)=0$ in $\C$, and in
particular, $\phi_0=\hbox{Re }\phi_0$ and
$\phi_k^1=\frac{N+1}{N+1+k}\hbox{Re }\phi_k$,
$\phi_k^2=\frac{N+1}{N+1+k}\hbox{Im }\phi_k$ are real solutions
for $k\geq 1$. We want to show that every solution $\phi$ of
$L(\phi)=0$ is a linear combination of $\phi_0$ and $\phi_k^i$,
$k\geq 1$ and $i=1,2$:
\begin{equation} \label{representation}
\phi=a_0 \phi_0+\sum_{k\geq 1}\left(a_k \phi_k^1+b_k^1
\phi_k^2\right).
\end{equation} The key idea is that, for $\rho$ small, the
functions $\phi_0(\rho e^{i\theta})$ and
$\frac{1}{\rho^k}\phi_k^i(\rho e^{i\theta)}$, $k\geq 1$ and
$i=1,2$, are very close to the real Fourier basis $1$,
$\cos(k\theta)$ and $\sin(k\theta)$ with $k\geq 1$, and then form
a complete set in $L^2(\partial B_\rho(0))$. Indeed, by a standard
integration by parts we can compute for $k\geq1$
$$\int_{S^1}\phi_0(\rho e^{i\theta})d\theta=1+O(\rho^{N+1})\:,\:\:\: \int_{S^1} \cos(k\theta)
\phi_0(\rho
e^{i\theta})d\theta=O(\frac{\rho^{N+1}}{k+1})\:,\:\:\:\int_{S^1}\sin(k\theta)\phi_0(\rho
e^{i\theta})d\theta=O(\frac{\rho^{N+1}}{k+1})$$ and
$$\int_{S^1}\cos(k\theta) \phi_j^1(\rho e^{i\theta})d\theta=\pi
\delta_{kj}+O(\frac{\rho^{N+1}}{(k+1)(j+1)})\:,\:\:\:\int_{S^1}\sin(k\theta)\phi_j^1(\rho
e^{i\theta})d\theta=O(\frac{\rho^{N+1}}{(k+1)(j+1)})$$ and
$$\int_{S^1} \cos(k\theta)
\phi_j^2(\rho e^{i\theta})d\theta=O(\frac{\rho^{N+1}}{(k+1)
(j+1)})\:,\:\:\:\int_{S^1} \sin(k\theta) \phi_j^2(\rho
e^{i\theta})d\theta=\pi
\delta_{kj}+O(\frac{\rho^{N+1}}{(k+1)(j+1)}).$$ Letting $\psi \in
L^2(\partial B_\rho(0))$ in the form
$$\psi(\rho e^{i\theta})=c_0+\sum_{k\geq 1}\left(c_k
\cos(k\theta)+d_k\sin(k\theta)\right),$$ we can compute
$$\tilde c_0:=\int_{S^1}\psi(\rho e^{i\theta}) \phi_0(\rho e^{i\theta})d\theta=c_0+
\rho^{N+1}O\left(\sum_{k\geq 0}\frac{|c_k|}{k+1}+\sum_{k\geq
1}\frac{|d_k|}{k+1}\right)$$ and
$$\tilde c_j:=\int_{S^1}\psi(\rho e^{i\theta}) \phi_j^1(\rho e^{i\theta})d\theta=c_j+
\frac{\rho^{N+1}}{j+1}O\left(\sum_{k\geq
0}\frac{|c_k|}{k+1}+\sum_{k\geq 1}\frac{|d_k|}{k+1}\right)$$
$$\tilde d_j:=\int_{S^1}\psi(\rho e^{i\theta}) \phi_j^2(\rho e^{i\theta})d\theta=d_j+
\frac{\rho^{N+1}}{j+1}\left(\sum_{k\geq
0}\frac{|c_k|}{k+1}+\sum_{k\geq 1}\frac{|d_k|}{k+1}\right).$$ We
consider the operator
$$T:(c_0,c_1,d_1,\dots)\in l_2 \:\to\: (\tilde c_0,\tilde c_1,\tilde,d_1,\dots) \in l^2.$$
We have shown so far that
$$\|T-\hbox{Id}\|\leq C\rho^{N+1} \left(\sum_{j\geq
0}\frac{1}{(j+1)^2}\right).$$ In conclusion, for $\rho$ small we
have that $T$ is an invertible operator, and by injectivity we
then deduce that, if $\psi \in L^2(\partial B_\rho(0))$ is so that
$$\int_{S^1}\psi(\rho e^{i\theta})\phi_0(\rho e^{i\theta})d\theta=
\int_{S^1}\psi(\rho e^{i\theta})\phi_k^j(\rho
e^{i\theta})d\theta=0 \quad \forall:k\geq 1,\:j=1,2,$$ then its
Fourier coefficients $c_j$ vanish and $\psi=0$. This means that,
for $\rho$ small, the space $L^2(\partial B_\rho(0))$ coincides
with the closure in $L^2-$norm of
$$\hbox{Span }\{\phi_0,\:\phi_k^j:\:k\geq 1,\:j=1,2\}.$$
In particular, every solution $\phi\in L^\infty(\C)$ can be
written on $\partial B_\rho(0)$, $\rho$ small, as
$$\phi(\rho e^{i\theta})=a_0 \phi_0(\rho e^{i\theta})+\sum_{k\geq
1}\left( a_k \phi_k^1(\rho e^{i\theta})+b_k \phi_k^2(\rho
e^{i\theta})\right),$$ for suitable $a_j$ and $b_j$. By regularity
theory $\phi \in C^\infty(\C)$, and then $\phi \mid_{\partial
B_\rho(0)} \in C^\infty(\partial B_\rho(0))$. Arguing as for the
Fourier coefficients, it is easily seen that $a_k$ and $b_k$ tend
to zero as $k \to +\infty$ faster than any power of $k$. In
particular, the function
$$\hat \phi(z)=a_0 \phi_0(z)+\sum_{k\geq 1}[a_k \phi_k^1(z)+b_k
\phi_k^2(z)]$$ is well defined, is in $C^\infty(\C)$ and satisfies
$L(\hat \phi)=0$ in $\C$. Since $\phi=\hat \phi$ on $\partial
B_\rho(0)$ and $\delta=\phi-\hat \phi$ satisfies $L(\delta)=0$ in
$\C$, an integration by parts yields to
$$\int_{B_\rho(0)}|\nabla
\delta|^2=\int_{B_\rho(0)}V\delta^2-\int_{B_\rho(0)}L(\delta)\delta=
\int_{B_\rho(0)}V\delta^2\leq C \rho^{2N}
\int_{B_\rho(0)}\delta^2,$$ where
$V(z)=\frac{8(N+1)^2|z|^{2N}}{(1+|z^{N+1}-c|^2)^2}.$ As soon as
$C\rho^{2N}<\lambda_1(B_\rho(0))$ ($\lambda_1$ being the first
eigenvalue of $-\Delta$ with Dirichlet boundary conditions), we
get that necessarily $\delta=0$ in $B_\rho(0)=0$. Then, for $\rho$
small we get that $\delta=0$ in $B_\rho(0)$, and by the strong
maximum principle $\delta=0$ in $\C$. So we have shown that
$$\phi(z)=a_0 \phi_0(z)+\sum_{k\geq 1}[a_k \phi_k^1(z)+b_k
\phi_k^2(z)]$$ in $\C$. Let us look now at the behavior of
$\phi(z)$ as $|z|\to +\infty$. Since the only bounded components
in $\phi(z)$ are $\phi_0$ and $\phi_{N+1}^1$, $\phi_{N+1}^2$, we
expect that $a_k=b_k=0$ for $k\not=0,N+1$. Also in this case we
will use that the components of $\phi$ are very close to the
Fourier basis as $|z|\to +\infty$.\\
Indeed, observe that
$$\frac{z^{N+1}\overline{z^{N+1}-c}}{1+|z^{N+1}-c|^2}=1+O(\frac{1}{|z|^{N+1}})\qquad\hbox{as
}|z|\to +\infty,$$ and then
$$\phi_k(z)=z^k\left(\frac{k-N-1}{N+1}+O(\frac{1}{|z|^{N+1}})\right)$$
at infinity. So we have that
\begin{eqnarray*}
&&\phi_0(z)=-1+O(\frac{1}{|z|^{N+1}})\:,\:\:\:
\phi_k^(z)=\frac{k-N-1}{N+1+k}|z|^k
\cos(k\theta)(1+O(\frac{1}{|z|^{N+1}})\\
&&\phi_k^2(z)=\frac{k-N-1}{N+1+k}|z|^k \sin(k\theta)
(1+O(\frac{1}{k|z|^{N+1}})). \end{eqnarray*} Let us now compute by
Cauchy-Schwartz inequality
\begin{eqnarray*}  \frac{1}{R}\int_{\partial B_R}|\phi|^2&=&\pi \left(\sum_{k\geq
0} \frac{(k-N-1)^2}{(N+1+k)^2} R^{2k}|a_k|^2+\sum_{k\geq 1}
\frac{(k-N-1)^2}{(N+1+k)^2} R^{2k}|b_k|^2\right)\\
&&+o\left(\sum_{k,j}R^{k+j}
(|a_k||a_j|+|b_k||b_j|+|a_k||b_j|)\right)\\
&&=\pi (1+o(1))\left(\sum_{k\geq 0} \frac{(k-N-1)^2}{(N+1+k)^2}
R^{2k}|a_k|^2+\sum_{k\geq 1} \frac{(k-N-1)^2}{(N+1+k)^2}
R^{2k}|b_k|^2\right)
\end{eqnarray*}
as $R\to +\infty$. Since $\phi \in L^\infty(\C)$, we have that
$\frac{1}{R}\int_{\partial B_R}|\phi|^2$ is bounded in $R$, and
then
$$\sum_{k\geq 0} \frac{(k-N-1)^2}{(N+1+k)^2}
R^{2k}|a_k|^2+\sum_{k\geq 1} \frac{(k-N-1)^2}{(N+1+k)^2}
R^{2k}|b_k|^2$$ is bounded in $R$. Then $a_k=0$ and $b_k=0$ for
$k\geq 1$ unless $k=N+1$. For a solution $\phi \in L^\infty(\C)$
of $L(\phi)=0$ we have then shown that
$$\phi(z)=a_0
\phi_0(z)+a_{N+1}\phi_{N+1}^1(z)+b_{N+1}\phi_{N+1}^2(z).$$ To
conclude, we need simply to rewrite $\phi_0$ and $\phi_{N+1}$ in a
more expressive way. We have that
$$\phi_0(z)=1-2\frac{z^{N+1}\overline{z^{N+1}-c}}{1+
|z^{N+1}-c|^2}=\frac{1-|z^{N+1}-c|^2}{1+
|z^{N+1}-c|^2}-2c\frac{\overline{z^{N+1}-c}}{1+ |z^{N+1}-c|^2}$$
and
$$\phi_{N+1}(z)=2 z^{N+1}\left(1-\frac{z^{N+1}\overline{z^{N+1}-c}}{1+ |z^{N+1}-c|^2}
\right)=2c \frac{1-|z^{N+1}-c|^2}{1+ |z^{N+1}-c|^2}+2
\frac{z^{N+1}-c}{1+ |z^{N+1}-c|^2}-2c^2
\frac{\overline{z^{N+1}-c}}{1+ |z^{N+1}-c|^2}.$$ In real form we
can then write that
$$\phi_0(z)=\frac{1-|z^{N+1}-c|^2}{1+
|z^{N+1}-c|^2}-2c_1 \hbox{Re }\frac{z^{N+1}-c}{1+ |z^{N+1}-c|^2}
-2c_2 \hbox{Im }\frac{z^{N+1}-c}{1+ |z^{N+1}-c|^2}$$ and
$$\phi_{N+1}^1(z)=c_1 \frac{1-|z^{N+1}-c|^2}{1+ |z^{N+1}-c|^2}+ (1-c_1^2+c_2^2)\hbox{Re }
\frac{z^{N+1}-c}{1+ |z^{N+1}-c|^2}-2c_1 c_2 \hbox{Im }
\frac{z^{N+1}-c}{1+ |z^{N+1}-c|^2}
$$
and
$$\phi_{N+1}^2(z)= c_2 \frac{1-|z^{N+1}-c|^2}{1+ |z^{N+1}-c|^2}-2c_1 c_2 \hbox{Re } \frac{z^{N+1}-c}{1+
|z^{N+1}-c|^2}+ (1+c_1^2-c_2^2) \hbox{Im }\frac{z^{N+1}-c}{1+
|z^{N+1}-c|^2},$$ where $c=c_1+ic_2$. In conclusion, the function
$\phi$ can be written as \begin{eqnarray*} \phi&=&(a_0
+a_{N+1}c_1+b_{N+1}c_2 )\frac{1-|z^{N+1}-c|^2}{1+
|z^{N+1}-c|^2}+[-2a_0c_1
+a_{N+1}(1-c_1^2+c_2^2)-2b_{N+1}c_1c_2]\hbox{Re
}\frac{z^{N+1}-c}{1+ |z^{N+1}-c|^2}\\
&&+[-2a_0 c_2-2a_{N+1}c_1c_2+b_{N+1}(1+c_1^2-c_2^2)] \hbox{Im
}\frac{z^{N+1}-c}{1+ |z^{N+1}-c|^2}
\end{eqnarray*}
and the Theorem is established.\qed


\begin{thebibliography}{99}

\bibitem{BaPa}S. Baraket, F. Pacard, Construction of singular
limits
for a semilinear elliptic equation in dimension
$2$, {\em Calc. Var. P.D.E.} {\bf 6} (1998), 1-38.

\end{thebibliography}
\end{document}